\documentclass[12pt]{article}
\usepackage[russian]{babel}
\usepackage{amsmath}
\usepackage[cp1251]{inputenc}
\usepackage{amscd}
\usepackage{amssymb}
\begin{document}

\noindent{\large\bf Представление типа Каратеодори с единичными
весами и связанные с ним аппроксимационные задачи}\footnote{This
work was supported by RFBR project 18-31-00312 mol\underline{\ \
}a.}\\

{\bf М. А. Комаров}\footnote{Department of Functional Analysis and
Its Applications, Vladimir State University, Gor$'$kogo str. 87,
600000 Vladimir, Russia\\ e-mail: kami9@yandex.ru}\\

{\bf Аннотация.} Произвольные $n$ комплексных чисел $a_{\nu-1}$,
$\nu=1,\dots,n$ при достаточно большом $n$ представляются в виде
степенных сумм:
$a_{\nu-1}=\lambda_1^\nu+\dots+\lambda_{2n+1}^\nu$, где
$\lambda_k$ --- попарно различные точки, лежащие на единичной
окружности. Рассматриваются приложения к задачам аппроксимации
суммами экспонент, $h$-суммами, к задаче выделения гармоник
сигнала. Результат базируется на оценке скорости аппроксимации
ограниченных аналитических в круге $|z|<1$ функций на замкнутых
подмножествах этого круга логарифмическими производными полиномов,
все нули которых лежат на окружности $|z|=1$. Наш результат это
модификация классического представления Каратеодори
$a_{\nu-1}=\sum_{k=1}^{n} X_k \lambda_k^\nu$, $\nu=1,2,\dots,n$
($X_k\ge 0$, $\lambda_k$ ---
попарно различные точки на единичной окружности).\\

MSC: 41A29, 41A21, 41A25, 41A20, 41A30\\

\noindent{\bf Введение}\bigskip

Классическая теорема Каратеодори устанавливает однозначную
разрешимость ``конечной'' тригонометрической проблемы моментов, а
именно \cite[\S 4.1]{Gren-Sz}: \medskip

{\it Пусть $a_0,a_1,\dots,a_{n-1}$ --- заданные комплексные
константы, не все равные нулю. Найдутся целое $m$, $1\le m\le n$,
и определённые постоянные $X_k,\lambda_k;\,k=1,2,\dots,m$, такие
что $X_k>0$, $|\lambda_k|=1$, $\lambda_k\ne \lambda_l$ при $k\ne
l$, и
\[a_{\nu-1}=\sum_{k=1}^m X_k \lambda_k^\nu, \qquad \eqno \nu=1,2,\dots,n.\]
Целое $m$ и константы $X_k,\lambda_k$ однозначно определены.}

\medskip

Это представление $n$ чисел $a_\nu$ содержит $\le 2n+1$ параметров
$X_k,\lambda_k,m$. В данной работе рассматриваем сходное
представление (теоремы 1.1, 1.2), но с {\it единичными весами}
$X_k=1$ и фиксированным $m=m_1=2n+1$:
\[a_{\nu-1}=\sum_{k=1}^{2n+1} \lambda_k^\nu, \eqno \nu=1,2,\dots,n\]
($|\lambda_k|=1$, $\lambda_k\ne \lambda_l$ при $k\ne l$). Здесь
ровно $2n+1$ параметров (при этом, $n$ не произвольно, например,
очевидно, $2n+1\ge \max\{|a_\nu|\}_{\nu=0}^{n-1}$.). Такие
представления\footnote{Мы не исследуем единственность
представления и минимизацию числа слагаемых в нём.} (степенными
суммами) естественным образом возникают в некоторых задачах
аппроксимации, например, в хорошо известной
\cite{Korev,Thompson,Rub-Suff} задаче аппроксимации
логарифмическими производными алгебраических полиномов, все корни
которых лежат на окружности (теорема 1.3), аппроксимации суммами
экспонент \cite{Beyl-Monz} и $h$-суммами \cite{D2008,ChD} (теоремы
2.1 и 2.2), в задаче о выделении
гармоник из тригонометрического многочлена \cite{DVas} (теорема 2.3).\\

\noindent{\bf 1. Основные результаты}\bigskip

{\bf 1.1.} Пусть $C$ это единичная окружность $|z|=1$, $D$
--- единичный круг $|z|<1$ и $\overline{D}=D\cup C$.
Пусть функция $f(z)$ аналитична в $D$. Обозначаем через
$s_n(z)=s_n(f;z)$ $n$-й полином Тейлора функции
$\exp\left(\int_0^z f(\zeta)d\zeta\right)$ и полагаем
\begin{equation}\label{P_n}
    P_n(z)=P_n(f;z):=s_n(z)+z^{N}\cdot \overline{s}_n(1/z), \qquad N:=2n+1.
\end{equation}
Здесь $\deg s_n(z)\le n$ и $\deg P_n(z)=N$, ибо $s_n(0)=1$. Мы
будем писать $f\in H^{\infty}(D)$, если $f$ --- ограниченная
аналитическая в $D$ функция. Степенные суммы $N$ точек
$\lambda_1,\dots,\lambda_N$, принадлежащих окружности $C$, будем
обозначать $S_\nu(\lambda)$:
\[S_\nu(\lambda)=S_\nu(\lambda_1,\dots,\lambda_N):=\lambda_1^{\nu}+\dots+\lambda_N^{\nu},
\qquad |\lambda_k|=1, \qquad \nu\in \mathbb{N}.\] Определим
$n_0(f)$ как наименьшее из натуральных $n$, таких что полином
$s_n(f;z)$ не имеет корней в $\overline{D}$.

\medskip

{\bf Теорема 1.1.} {\it Пусть $f\in H^{\infty}(D)$,
$f(z):=\sum_0^\infty (-a_j) z^j$, $n_0=n_0(f)$. Существует номер
$n_1=n_1(f)\ge n_0$ такой, что при $n\ge n_1:$

$1)$ все $N$ корней $z_k$ полинома $P_n(f;z)$ $($см.
$(\ref{P_n}))$ простые и принадлежат $C$,

$2)$ для попарно различных точек
\[\lambda_k:=z_k^{-1}, \qquad |\lambda_k|=1, \eqno k=1,\dots,N,\]
принадлежащих единичной окружности, верны соотношения
\[S_{j+1}(\lambda)-a_{j}=0, \eqno j=0,\dots,n-1;\]
\[|S_{j+1}(\lambda)-a_j|<\frac{r^{n+1}(r-\varepsilon)^{-j}}
{2\varepsilon(1-r)}, \eqno j\ge n,\] какими бы ни были числа $r$ и
$\varepsilon$, $0<\varepsilon<r<1$. При этом, для выполнения
равенств $S_{j+1}(\lambda)-a_{j}=0$, $j=0,\dots,n-1$ достаточно,
чтобы $n\ge n_0$.}

\bigskip

Обозначим $\lambda(n;\{a_j\})$ набор
$\lambda_1,\dots,\lambda_{2n+1}$, определённый в теореме 1.1, 2).

\bigskip

{\bf Теорема 1.2.} {\it Если $\lambda=\lambda(n;\{a_j\})$, где
$n\ge 1$ и все $|a_j|\le (j+2)^{-2}$, то
\[S_{j+1}(\lambda)-a_{j}=0, \eqno j=0,\dots,n-1,\]
\[|S_{j+1}(\lambda)-a_j|\le \frac{r^{n-j}}{1-r^{n+1}}\left(15n+\frac{30}{1-r}\right), \eqno
j\ge n,\] где $r\in (0,1)$ любое.}

\medskip

{\bf Замечание 1.1.} В \cite{Gavrilov} дан некоторый рекуррентный
алгоритм подбора конечного числа точек $\lambda_k$, таких что
$|\lambda_k|=1$ и $\sum_k \lambda_k^m=a_{m-1}$, $1\le m \le n$,
для произвольно заданных чисел $a_{\nu}$, но не устанавливается
достаточное количество этих точек, $N$.

\medskip

{\bf 1.2.} Теоремы 1.1, 1.2 следуют из теоремы 1.3,
устанавливающей скорость равномерной аппроксимации ограниченных
аналитических в $D$ функций суммами вида
\begin{equation}\label{spf_N}
    \sum_{k=1}^{n}(z-z_{k})^{-1}, \qquad |z_{k}|=1,
\end{equation}
т.е. логарифмическими производными $C$-{\it полиномов} (так
принято называть полиномы, все нули которых лежат на $C$). Задаче
аппроксимации рациональными функциями вида (\ref{spf_N}) с
ограничениями на полюсы посвящено немало работ (см.
\cite{Korev,Thompson,Rub-Suff,Chui-Shen}). Например, в
\cite{Thompson,Rub-Suff} показано, что если $f\in H^\infty(D)$, то
существует последовательность рациональных функций вида
(\ref{spf_N}), сходящаяся к $f(z)$ равномерно на замкнутых
подмножествах $D$. Однако, оценок скорости такой аппроксимации до
настоящего времени не было. Теорема 1.3 доказана в \S 4.

\medskip

{\bf Теорема 1.3.} {\it Пусть $f\in H^\infty(D)$, $n\ge n_0(f)$ и
$N:=2n+1$. Тогда $P_n$ $($см. $(\ref{P_n}))$ это $C$-полином с
попарно различными корнями и
\begin{equation}\label{interp}
    P_n'(z)/P_n(z)-f(z)=O(z^n), \qquad z\to 0.
\end{equation}
Более того, найдётся номер $n_1=n_1(f)$, такой что если $n\ge
n_1$, то в любом круге $|z|\le a<1$ при каждом $\varepsilon\in
(0,1-a)$
\begin{equation}\label{result}
    \left|\frac{P_n'(z)}{P_n(z)}-f(z)\right|<
    \frac{(a+\varepsilon)^{n+1}}{2\varepsilon(1-a-\varepsilon)}.
\end{equation}
Если для $f(z)=\sum_{j\ge 0} f_j z^j$ выполнено $|f_j|\le
(j+2)^{-2}$, $j=0,1,2,\dots$, то при каждом $n\ge 1$ имеем более
сильную оценку}
\begin{equation}\label{result(|fj|<...)}
    \left|\frac{P_n'(z)}{P_n(z)}-f(z)\right|\le
    \frac{15|z|^n}{1-|z|^{n+1}} \left(n+\frac{2}{1-|z|}\right),
    \qquad |z|<1.
\end{equation}

{\it Доказательство теоремы 1.1.} Положим $n\ge n_0(f)$. Разложим
логарифмическую производную $P_n'/P_n$ полинома (\ref{P_n}) в ряд
по степеням $z$:
\[\frac{P_n'(z)}{P_n(z)}=\sum_{k=1}^N \frac{1}{z-z_k}=
-\sum_{k=1}^N \frac{z_k^{-1}}{1-z_k^{-1}z}=-\sum_{j=0}^{\infty}
z^j\sum_{k=1}^N z_k^{-j-1}, \quad |z_k|=1, \quad |z|<1.\] Отсюда
при $\lambda_k:=z_k^{-1}$, $S_j=S_j(\lambda)$ и
$f(z)=-\sum_0^{\infty} a_j z^j$ получаем
\[\Delta(z):=\frac{P_n'(z)}{P_n(z)}-f(z)\equiv -\sum_{j=0}^{\infty}(S_{j+1}-a_j)z^j, \qquad |z|<1.\]
Ввиду $\Delta(z)=O(z^n)$ (см. (\ref{interp})) имеем
$S_{j+1}-a_j=0$, $j=0,\dots,n-1$.

При $n\ge n_1(f)$ верна оценка (\ref{result}). Величины
$S_{j+1}-a_j$, $j\ge n$, оцениваются с помощью неравенств Коши как
$j$-й коэффициент Тейлора функции $-\Delta$, аналитической в круге
$|z|\le r-\varepsilon$, $0<\varepsilon<r<1$:
\[|S_{j+1}-a_j|\le (r-\varepsilon)^{-j}\max_{|z|=r-\varepsilon}|\Delta(z)|<
\frac{r^{n+1}(r-\varepsilon)^{-j}}{2\varepsilon(1-r)}, \qquad j\ge
n.\]

{\it Доказательство теоремы 1.2.} При $|a_j|\le (j+2)^{-2}$
согласно (\ref{result(|fj|<...)}) и оценкам Коши при каждых $n\ge
1$ и $r\in (0,1)$ имеем $S_{j+1}-a_j=0$, $j=0,\dots,n-1$ и
\[|S_{j+1}-a_j|\le r^{-j}\max_{|z|=r}|\Delta(z)|\le
\frac{r^{n-j}}{1-r^{n+1}}\left(15n+\frac{30}{1-r}\right), \qquad
j\ge n.\]

{\bf Замечание 1.2.} В \cite{DD,Kos,K-IzvRAN-2017} изучалась
аппроксимация суммами $\sum_{k=1}^n(z-z_k)^{-1}$ со {\it
свободными} полюсами (см. также библиографию в
\cite{K-IzvRAN-2017}). Скорость (\ref{result}) вообще говоря
сравнима со скоростью $d_n(f,K)$ приближения функций $f$ в круге
$K:|z|\le a$ логарифмическими производными $Q'(z)/Q(z)$ со
свободными полюсами (здесь $\deg Q\le n$):
\begin{equation}\label{dn<a^n}
    \limsup_{n\to \infty} \sqrt[n]{d_n(f,K)}\le a.
\end{equation}
Оценка (\ref{dn<a^n}) вытекает из \cite{Kos} (см. также
\cite{DD}), где для задачи аппроксимации такими логарифмическими
производными ({\it наипростейшими дробями}) построен аналог
теоремы Уолша.\\

\noindent{\bf 2. Приложения к задачам аппроксимации}
\bigskip

{\bf 2.1.} {\it Аппроксимация суммами комплексных экспонент.}

\medskip

{\bf Теорема 2.1.} {\it Пусть \[f(z):=\sum_{0}^\infty p_j
\frac{z^j}{j!}, \qquad |p_j|\le (j+2)^{-2}, \qquad N=2n+1, \quad
n\in \mathbb{N}.\] При каждом $n\ge 1$ найдутся попарно различные
числа $\lambda_1,\dots,\lambda_N$, $|\lambda_k|=1$, такие что при
любом $r\in(0,1)$ во всей плоскости $\mathbb{C}$ имеем}
\[\left|\sum_{k=1}^N \lambda_k e^{\lambda_k z}-f(z)\right|\le
\frac{|z|^n}{n!}\frac{15}{1-r^{n+1}}
\left(n+\frac{2}{1-r}\right)\left(1+\frac{|z|e^{|z|/r}}{rn+r}\right),
\qquad z\in \mathbb{C}.\]

{\it Доказательство.} Рассмотрим разложение
\[\Delta(z):=\sum_{k=1}^N \lambda_k e^{\lambda_k z}-f(z)=\sum_{j=0}^\infty
(S_{j+1}-p_j)\frac{z^j}{j!}, \qquad S_{j+1}=S_{j+1}(\lambda).\]
Положим $n\ge 1$, $\lambda=\lambda(n;\{p_j\})$. Поскольку все
$|p_j|\le (j+2)^{-2}$, по теореме 1.2 имеем при любом $0<r<1$
\[|\Delta(z)|\equiv \left|\sum_{j=n}^\infty
(S_{j+1}-a_j)\frac{z^j}{j!}\right|\le
\frac{|z|^n}{n!}\frac{15}{1-r^{n+1}}
\left(n+\frac{2}{1-r}\right)\sum_{j=n}^\infty
\frac{|z|^{j-n}}{r^{j-n}}\frac{n!}{j!}.\] Последняя сумма не
превосходит $1+|z|e^{|z|/r}/(rn+r)$ при всех $z$. Теорема 2.1
доказана.

\smallskip

{\bf 2.2.} {\it Аппроксимация посредством $h$-сумм.} В
\cite{D2008} предложен метод аппроксимации так называемыми
$h$-{\it суммами}
\begin{equation}\label{H_n}
    H_n(\lambda;z)=\sum_{k=1}^n \lambda_k h(\lambda_k z), \qquad
    h(z):=\sum_{j=0}^\infty h_j z^j, \qquad 0<|h_j|\le M,
\end{equation}
где $h(z)$ --- фиксированная аналитическая в $D$ функция.
Аппроксимация производится за счёт выбора параметров $\lambda_k$.
$h$-суммы это естественное обобщение рациональных функций
(\ref{spf_N}) (где $h(z):=(z-1)^{-1}$, $\lambda_k:=z_k^{-1}$) и
сумм экспонент $\sum \lambda_k e^{\lambda_k z}$ ($h(z)=e^z$).

В \cite{D2008} показано, что скорость приближения суммами
(\ref{H_n}) {\it без ограничений на $\lambda_k$} геометрическая.
Здесь мы докажем аналогичный результат при $|\lambda_k|=1$.

\medskip

{\bf Теорема 2.2.} {\it Пусть $h$ --- фиксированная аналитическая
в $D$ функция $($см. $(\ref{H_n}))$,
\[f(z):=\sum_{0}^\infty f_j z^j, \quad
F(z):=\sum_{0}^\infty \frac{f_j}{h_j}z^j, \quad F\in H^\infty(D),
\quad N=2n+1, \quad n\in \mathbb{N}.\] При $n\ge n_1(-F)$ найдутся
попарно различные числа $\lambda_1,\dots,\lambda_N$,
$|\lambda_k|=1$, такие что
\[|H_N(\lambda;z)-f(z)|\le   |z|^{n}\frac{(5-|z|)^{n+1}}{4^{n-1}}
\frac{2M}{3}\frac{3+|z|}{(1-|z|)^{4}}, \qquad |z|<1;\] здесь
$(5|z|-|z|^2)/4<1$.}

\medskip

{\it Доказательство.} В условиях теоремы функция $f$
аналитична\footnote{В самом деле, $\overline{\lim}_{j\to \infty}
\sqrt[j]{|f_j/h_j|}\le 1$ и $|f_j/h_j|\ge |f_j|/M$ по условию,
поэтому $\overline{\lim}_{j\to \infty} \sqrt[j]{|f_j|}\le
\overline{\lim}_{j\to \infty} \sqrt[j]{M |f_j/h_j|}\le 1$.} в $D$.
Пусть $n\ge n_1(-F)$. Положим $\lambda=\lambda(n;\{f_j/h_j\})$.
Функция $H_N(\lambda;z)$ также аналитична в $D$. Имеем
\[\Delta(z):=H_N(\lambda;z)-f(z)\equiv \sum_{j=n}^\infty (h_j S_{j+1}-f_j)z^j, \qquad
S_{j+1}=S_{j+1}(\lambda),\]
\[|\Delta(z)|\le \frac{Mr^{n+1}}{1-r}\sum_{j=n}^\infty
\frac{|z|^j/(2\varepsilon)}{(r-\varepsilon)^j}=
\frac{M}{r-\varepsilon-|z|}\frac{r^{n+1}}{1-r}
\frac{|z|^n/(2\varepsilon)}{(r-\varepsilon)^{n-1}}, \quad
|z|<r-\varepsilon<r<1.\] Выберем $r=1-\delta^2$,
$\varepsilon=(1-\delta)\delta$, где $\delta=(1-|z|)/4$. Тогда
$1-r=\delta^2$, $r-\varepsilon=1-\delta$,
$r-\varepsilon-|z|=3\delta$. После простых преобразований получаем
указанную оценку. Теорема 2.2 доказана.

\smallskip

{\bf Замечание 2.1.} В \cite{ChD} изучается дальнейшее обобщение
аппаратов (\ref{spf_N}) и (\ref{H_n}) --- {\it
ампли\-туд\-но-частотные суммы}
\[\sum_{k=1}^n \mu_k h(\lambda_k z),\]
где $h$ такая же, как в теореме 2.2, а амплитуды $\mu_k$ не
зависят от $\lambda_k$. Эти аппараты применяются для численного
дифференцирования, интегрирования, экстраполяции (см. \cite{ChD} и
библиографию в \cite{ChD}). Рассмотрим специальные
амплитудно-частотные суммы, где все $\mu_k=1$, $|\lambda_k|=1$, а
$h_0=h(0)=0$:
\begin{equation}\label{H_n^(1)}
    H_n^{(1)}(\lambda;z)=\sum_{k=1}^n h(\lambda_k z), \qquad
    |\lambda_k|=1, \quad h(z):=\sum_{j=1}^\infty h_j z^{j},
    \quad 0<|h_j|\le M.
\end{equation}

{\bf Следствие 2.1.} {\it Пусть $h$ --- фиксированная
аналитическая в $D$ функция $($см. $(\ref{H_n^(1)}))$,
\[f(z):=\sum_{1}^\infty f_j z^j, \quad F_1(z):=\sum_{1}^\infty
\frac{f_j}{h_j}z^{j-1}, \quad F_1\in H^\infty(D), \quad N=2n+1,
\quad n\in \mathbb{N}.\] При $n\ge n_1(-F_1)$ найдутся попарно
различные числа $\lambda_1,\dots,\lambda_N$, $|\lambda_k|=1$,
такие что} \[|H_N^{(1)}(\lambda;z)-f(z)|\le
\frac{(5|z|-|z|^2)^{n+1}}{4^{n-1}}
\frac{2M}{3}\frac{3+|z|}{(1-|z|)^{4}} \qquad |z|<1.\]

В самом деле, имеем $h(z)\equiv z \tilde{h}(z)$, $f(z)\equiv z
\tilde{f}(z)$ и $H_N^{(1)}(\lambda;z)-f(z)\equiv
z[\tilde{H}_N(\lambda;z)-\tilde{f}(z)]$, где $\tilde{h},\tilde{f}$
аналитичны в $D$, а $\tilde{H}_N(\lambda;z):=\sum_{1}^N \lambda_k
\tilde{h}(\lambda_k z)$ это сумма вида (\ref{H_n})
($\tilde{h}$-сумма). Остаётся применить теорему 2.2 к разности
$\tilde{H}_N-\tilde{f}$.

\smallskip

{\bf 2.3.} {\it Выделение гармоник из тригонометрического
многочлена.} Практический интерес имеет задача выделения гармоник
из тригонометрического многочлена (сигнала) $T_n(t)=\sum_{m=-n}^n
c_m e^{itm}$, $n\ge 1$. Будем считать, что сигнал вещественный с
нулевым средним, т.е. $c_{-m}=\overline{c_m}$, $c_0=0$. Таким
образом,
\begin{equation}\label{T_n}
    T_n(t)\equiv 2{\rm Re}\sum_{m=1}^n c_m e^{itm}\equiv \sum_{m=1}^n
    \tau_m(t), \quad \tau_m(t):=a_m\cos mt+b_m\sin mt,
\end{equation}
\[a_m,b_m\in \mathbb{R}, \qquad a_m-ib_m=2c_m.\]
Недавно в \cite{DVas} поставлена и решена задача о выделении
отдельных гармоник $\tau_\nu$ с помощью {\it универсальных
вещественных амплитудно-фазовых операторов} вида
\[\tau_\nu(t)\equiv A_{N'}(t), \qquad A_{N'}(t):=\sum_{k=1}^{N'} X_k T_n(t-t_k),
\quad X_k,t_k\in \mathbb{R}.\] Универсальность оператора означает,
что его {\it параметры $X_k,t_k,N'$ зависят только от $n$ и $\nu$
--- номера выделяемой гармоники, --- и не
зависят от многочлена\footnote{В частности, оператор выделяет
гармоники в случае любых $a_m,b_m\in \mathbb{C}$.} $T_n$}.

В \cite{DVas} решение найдено в явном виде; показано, что порядок
оператора $N'\le n+1$, причём все $X_k$ одного знака и
$|X_1+\dots+X_{N'}|\le 2$. Этот подход гораздо проще известных
спектральных методов (см. библиографию в \cite{DVas}), ибо не
содержит операций интегрирования. Он позволил получить
\cite{DVas,DD2018} точные оценки гармоник и пар гармоник через
sup-норму сигнала $T_n$.

Мы построим аналогичные универсальные вещественные операторы
порядка $N=2n+1$, все амплитуды которых равны $X_k=1$, а фазы
$t_k=t_k(n,\nu)$ попарно различны. Предполагаем $n\ge 2$ (только
этот случай и представляет интерес).

\medskip

{\bf Теорема 2.3.} {\it Пусть целое $n\ge 2$,
$\nu\in\{1,2,\dots,n\}$ и $N=2n+1$. Пусть $z_k$ $(k=1,2,\dots,N)$
это корни многочлена
\begin{equation}\label{s_n(z^nu;z)}
    s_n(z)+z^N s_n(1/z), \qquad s_n(z):=\sum_{j=0}^q
    \frac{(-z^\nu)^{j}}{\nu^j j!}, \qquad q:=[n/\nu].
\end{equation}
Числа $z_k$ лежат на единичной окружности, числа $t_k:=\arg z_k\in
[0,2\pi)$ попарно различны, и для любого полинома $T_n(t)$ вида
$(\ref{T_n})$ справедливо тождество
\begin{equation}\label{tau_nu=...}
    \tau_\nu(t)\equiv \Theta_N(n,\nu;t), \qquad
    \Theta_N(n,\nu;t):=\sum_{k=1}^N T_n(t-t_k),
\end{equation}
где $\tau_\nu$ --- гармоника с номером $\nu$ полинома $T_n$.}

\medskip

{\it Доказательство.} В теореме 1.1 возьмём $f(z)=-z^{\nu-1}$
(т.е. $a_{\nu-1}=1$, $a_j=0$ при $j\ne \nu-1$). Очевидно, полином
$P_n(f;z)$ (см. $(\ref{P_n})$) имеет вид (\ref{s_n(z^nu;z)}).
Покажем, что
\[n_0(-1)=2; \qquad n_0(-z^{\nu-1})=1, \qquad \nu=2,3,\dots.\]
Действительно, для $f(z)=-1$ имеем $s_n(z)=\sum_0^n (-z)^j/j!$, и
только $s_1(z)$ имеет корень в $\overline{D}$, а для
$f(z)=-z^{\nu-1}$ ($\nu=2,3,\dots$) при любом $n\in \mathbb{N}$
полином $s_n(z)$ (см. (\ref{s_n(z^nu;z)})) не имеет корней в
$\overline{D}$: в случае $q=0$ это очевидно, а в случае $q\ge 1$
\[|s_n(z)|>1-\sum_{j=1}^\infty \frac{1}{\nu^j j!}>
1-\frac{e^{1/\nu}}{\nu}>0, \qquad |z|\le 1.\] Итак, $n_0(f)\le 2$,
следовательно, для любого $n\ge 2$ корни полинома $P_n(f;z)$ имеют
вид $z_k=e^{it_k}$, $t_k\in [0,2\pi)$, $t_k\ne t_j$ ($k\ne j$) и
выполняются равенства $S_\nu(\lambda)=1$, $S_j(\lambda)=0$,
$j\in\{1,2,\dots,n\}\backslash\{\nu\}$, где $\lambda$ --- набор
чисел $\lambda_k:=z_k^{-1}=e^{-it_k}$. Тем самым,
\[\sum_{k=1}^N T_n(t-t_k)=
2{\rm Re}\sum_{k=1}^N \sum_{j=1}^{n} c_j e^{i(t-t_k)j}=2{\rm
Re}\sum_{j=1}^{n}c_j e^{itj} S_j(\lambda)= 2{\rm Re}\ c_\nu e^{i
t\nu}\equiv \tau_\nu(t).\] Теорема 2.3 доказана.

\medskip

{\bf Следствие 2.2.} {\it Коэффициенты Фурье многочлена $T_n(t)$
вида $(\ref{T_n})$, $n\ge 2$, можно вычислить по формулам
\[a_\nu =\sum_{k=1}^{2n+1} T_n(-t_k), \qquad b_\nu
=\sum_{k=1}^{2n+1} T_n\left(\frac{\pi}{2\nu}-t_k\right), \eqno
\nu=1,\dots,n;\] числа $t_k$ определены в теореме $2.3$.}

\medskip

{\bf Замечание 2.2.} Очевидно, из цитированной во Введении теоремы
Каратеодори следует существование единственного
амплитудно-фазового оператора $A_n(t)$ (в смысле \cite{DVas})
порядка $\le n$, выделяющего заданную линейную комбинацию гармоник
полинома $T_n(t)\not\equiv 0$. В \cite{DD2018} найдены явные
формулы операторов, выделяющих суммы гармоник $\tau_1+\tau_n$ и
$\tau_2+\tau_n$. Нетрудно видеть, что задача выделения линейной
комбинации {\it начальных гармоник} (а именно эти гармоники и
несут основную информацию сигнала) решается и с помощью сумм вида
$\Theta_N$: {\it для любого $\tilde{n}\in \mathbb{N}$ и любых
чисел $\gamma_m\in \mathbb{C}$ $(m=1,2,\dots,\tilde{n})$
существует номер $n^*=n^*(\{\gamma_m\})\ge \tilde{n}$, такой что
при $n\ge n^*$ найдутся $N=2n+1$ попарно различных чисел $t_k\in
[0,2\pi)$, таких что \[\sum_{m=1}^{\tilde{n}}
\gamma_m\tau_m(t)\equiv \sum_{k=1}^N T_n(t-t_k),\] каким бы ни был
полином $T_n(t)$ вида $(\ref{T_n})$ степени $\le n$.}
Доказательство вполне аналогично доказательству теоремы 2.3; здесь
$n^*:=n_0(f)$, где
$f(z):=\sum_{j=1}^{\tilde{n}} (-\gamma_j) z^{j-1}$. \\

\noindent{\bf 3. Вспомогательные результаты}\bigskip

{\bf 3.1.} Свойства полиномов (\ref{P_n}) вытекают из следующей
леммы, полученной в \cite{Rub} (случай $m=0$ см. в
\cite[p.\,108]{Polya-Szego}).

\medskip

{\bf Лемма 3.1} \cite{Rub}. {\it Пусть
$Q(z)=q_0(z-u_1)\dots(z-u_q)$, $q_0\ne 0$, --- полином степени
$q\ge 1$, и
$Q^*(z):=z^q\overline{Q}(1/z)=\overline{q}_0(1-\overline{u}_1
z)\dots (1-\overline{u}_q z)$. Если $Q(z)$ не имеет нулей в
$\overline{D}$, то при любом $m=0,1,2,\dots$ сумма $Q(z)+z^m
Q^*(z)$ это $C$-полином, имеющий ровно $q+m$ попарно различных
корней, и $|Q^*(z)|\le |Q(z)|$ в $\overline{D}$.}

\medskip

Для доказательства леммы достаточно рассмотреть уравнение
\begin{equation}\label{Phi=...}
    \Phi(z)=-q_0/\overline{q}_0, \qquad
    \Phi(z):=\frac{q_0}{\overline{q}_0}\frac{z^m Q^*(z)}{Q(z)}\equiv
    z^m \prod_{j=1}^q \frac{1-\overline{u}_j z}{z-u_j}, \qquad
    |u_j|>1.
\end{equation}
Модуль каждого сомножителя в последнем произведении меньше, равен
либо больше 1, если и только если $|z|<1$, $|z|=1$ либо $|z|>1$,
соответственно, поэтому все корни уравнения (\ref{Phi=...}) (и
нули многочлена $Q(z)+z^m Q^*(z)$) лежат на $C$, а в круге
$\overline{D}$ имеем $|\Phi(z)|\le 1$ и $|Q^*(z)|\le |Q(z)|$.

Функция $\Phi(z)$ в $D$ имеет $q+m$ корней с учётом кратностей и
не имеет полюсов, следовательно, при обходе точкой $z$ единичной
окружности аргумент $\arg \Phi(z)$ получает приращение $2\pi
(q+m)$. Поэтому уравнение (\ref{Phi=...}) имеет не менее $q+m$
различных корней. Но $\deg (Q(z)+z^m Q^*(z))\le q+m$, значит,
корней ровно $q+m$.

{\bf Замечание 3.1.} Если $Q(z)\equiv a_0={\rm const}\ne 0$, то
$Q(z)+z^m Q^*(z)\equiv a_0+z^m \overline{a}_0$ это $C$-полином
только при $m>0$.

\medskip

{\bf 3.2.} При доказательстве теоремы 1.2 нам потребуется при
$\sigma=1$ следующая

\medskip

{\bf Лемма 3.2.} {\it Пусть $f(z)=\sum_{0}^{\infty} f_{k} z^k$ и
$|f_k|\le (k+2)^{-1-\sigma}$, $\sigma>0$. Тогда для тейлоровских
коэффициентов функции $g(z)=\exp\left(\int_0^z
f(\zeta)d\zeta\right)=1+\sum_{1}^{\infty} g_k z^k$ имеем}
\[|g_1|\le 2^{-1-\sigma}, \qquad |g_k|<(k+1)^{-1-\sigma}, \qquad k=2,3,\dots.\]

Для доказательства положим $f_k=(-1)^k \beta_{k+1}$. По условию,
$|\beta_k|\le (k+1)^{-1-\sigma}$, $k=1,2,\dots$. Воспользуемся
рекуррентными формулами \cite[лемма 1]{Ch}:
\[g_1=\beta_1, \qquad g_k=\frac{(-1)^{k+1}}{k}\left(\beta_k+
\sum_{j=1}^{k-1}(-1)^j \beta_{k-j}g_j \right), \quad
k=2,3,\dots.\] Отсюда $|g_1|\le 2^{-1-\sigma}$, $|g_2|\le
(|\beta_2|+|\beta_1||g_1|)/2\le
(3^{-1-\sigma}+4^{-1-\sigma})/2<3^{-1-\sigma}$. Пусть $k\ge 3$.
Если $|g_j|<(j+1)^{-1-\sigma}$ для $j=2,\dots,k-1$, то
\[k|g_k|\le |\beta_k|+\sum_{j=1}^{k-1} |\beta_{k-j}||g_j|<
\sum_{j=0}^{k-1} \frac{1}{(k+1-j)^{1+\sigma}(j+1)^{1+\sigma}}<
\frac{k}{(k+1)^{1+\sigma}}\] (воспользовались тем, что
$(k+1-j)(j+1)>k+1$ при всех $j$ от 1 до $k-1$). Лемма 3.2 доказана
по индукции.

\smallskip

{\bf 3.3.} Нам потребуются также хорошо известные формулы
\begin{equation}\label{sum_k^(-s)}
    \sum_{k=1}^{\infty} \frac{1}{k^{2}}=\frac{\pi^2}{6},
    \qquad \sum_{k=n+1}^{\infty} \frac{1}{k^{2}}\le \frac{1}{n} \qquad (n\in \mathbb{N}).
\end{equation}

\noindent{\bf 4. Доказательство теоремы 1.3}\bigskip

{\bf 4.1.} {\it Случай $f\in H^{\infty}(D)$.} Положим
$g(z)=\exp\left(\int_0^z f(\zeta)d\zeta\right)$,
\[g(z)=s_n(z)+R_{n}(z), \qquad s_n(z)=1+\sum_{k=1}^n g_k z^k, \qquad
R_n(z)=\sum_{k=n+1}^{\infty} g_k z^k, \qquad n\ge 1.\]

Производная $g'=gf$ ограничена и аналитична в $D$. В частности,
$g'(z)$ принадлежит классу Харди $H^1(D)$ и ряд $\sum |g_k|$
сходится \cite{Hardy}. Отсюда следует, что существует номер
$n_0=n_0(f)<\infty$, определённый в \S 1.1. В самом деле, функция
$\varphi(x,y)={\rm Re} \int_0^z f(\zeta)d\zeta$ ограничена в $D$,
поэтому $M_0:=\inf_{D} |g(z)|=\inf_{D} \exp \varphi(x,y)>0$. Ясно,
что $n_0\le n'$, где номер $n'$ такой, что $\sum_{n'+1}^{\infty}
|g_k|\le M_0/2$.

Итак, полиномы $s_n(z)\equiv g(z)-R_n(z)$, $n\ge n_0$, не имеют
нулей в $\overline{D}$.

При $n\ge n_0$ положим $p(z)=s_n^*(z)=z^{q}\overline{s}_n(1/z)$,
где $q=\deg s_n(z)$, $0\le q\le n$. По лемме 3.1 и замечанию 3.1
имеем $|p(z)|\le |s_n(z)|$ в $\overline{D}$, а сумма
\[P(z)=P_n(z):=s_n(z)+z^{m}p(z), \qquad m:=2n+1-q,\]
это $C$-полином ($\deg P_n=2n+1$, ибо $s_n(0)=1$ и $\deg p(z)=q$).

Перепишем $P$ в виде $P(z)\equiv g(z)+z^{m}p(z)-R_n(z)$. Имеем
\[P'(z)=g(z)f(z)+m z^{m-1}p(z)+z^{m}p'(z)-R_n'(z),\]
\begin{equation}\label{P'/P-f}
    \frac{P'(z)}{P(z)}-f(z)= \frac{z^{m-1}(m
    -zf(z))p(z)+z^{m}p'(z)-R_n'(z)+f(z)R_n(z)}{P(z)}.
\end{equation}
Из (\ref{P'/P-f}) следует, что равенство (\ref{interp}) выполнено
при всех $n\ge n_0$, ибо $P(0)=1$, $m\ge n+1$ и
$R_n(z)=O(z^{n+1})$. Далее полагаем $n\ge n'$.

Пусть $M_1=\sum_0^\infty |g_k|$ ($M_1<\infty$). Согласно
неравенству Бернштейна,
\begin{equation}\label{|p'|}
    |p'(z)|\le (\deg p)\max_{|\zeta|=1}|p(\zeta)|
    \le q\max_{|\zeta|=1}|s_n(\zeta)|\le n M_1, \qquad |z|\le 1.
\end{equation}
Далее, поскольку $|g_k|\le M_0/2$ \ ($k\ge n'+1$), имеем
\begin{equation}\label{|pn|,|Rn|}
    |R_n(z)|\le \frac{M_0}{2}\frac{|z|^{n+1}}{1-|z|},
\end{equation}
\begin{equation}\label{|P|>}
    |P(z)|\ge |g(z)|-|z^{m}p(z)-R_n(z)|\ge
    M_0-n M_1|z|^m-\frac{M_0}{2}\frac{|z|^{n+1}}{1-|z|}
\end{equation}
при $|z|<1$. Если $|z|<r<1$ и функция $F$ аналитична в $D$, то
\[2\pi|F'(z)|=\left|\int_{|\zeta|=r}\frac{F(\zeta)d\zeta}{(\zeta-z)^2}\right|
\le \max_{|\zeta|=r}|F(\zeta)|
\int_{|\zeta|=r}\frac{|d\zeta|}{|\zeta-z|^2}=
\max_{|\zeta|=r}|F(\zeta)|\frac{2\pi r}{r^2-|z|^2}\] (мы применили
интеграл Пуассона), и если $r=a+\varepsilon<1$
($a,\varepsilon>0$), то
\[|F'(z)|\le \frac{\tau}{\varepsilon}\cdot\max_{|\zeta|=a+\varepsilon}|F(\zeta)|,
\qquad |z|\le a, \qquad
\tau:=\frac{a+\varepsilon}{2a+\varepsilon}<1.\] Отсюда ввиду
(\ref{|pn|,|Rn|}) имеем при любом $a\in (0,1)$ и $\varepsilon\in
(0,1-a)$
\begin{equation}\label{|p'|<}
    |R_n'(z)|\le \frac{\tau M_0(a+\varepsilon)^{n+1}}{2\varepsilon(1-a-\varepsilon)},
    \qquad |z|\le a.
\end{equation}
Из (\ref{P'/P-f})---(\ref{|p'|<}) видно, что при $|z|\le a$ и
$n\ge n'$
\[\left|\frac{P'(z)}{P(z)}-f(z)\right|\le
\frac{\tau(a+\varepsilon)^{n+1}}{2\varepsilon(1-a-\varepsilon)}\cdot\psi_n,
\qquad \psi_n=\psi_n(a,\varepsilon,M_0,M_1),\] где $\psi_n\to 1$
при $n\to \infty$. Поскольку $\tau<1$, найдётся номер $n_1\ge n'$,
такой что $\tau\psi_n<1$ при всех $n\ge n_1$. Получили оценку
(\ref{result}).

\medskip

{\bf 4.2.} {\it Случай} $|f_k|\le (k+2)^{-2}$, $k=0,1,2,\dots$. В
этом случае, по лемме 3.2 коэффициенты Тейлора функции $g(z)$
удовлетворяют оценке $|g_k|\le (k+1)^{-2}$. В частности (см.
(\ref{sum_k^(-s)})), при любом $n\ge 1$
\[|g_1|+\dots+|g_n|\le \sum_2^{n+1} \frac{1}{k^2}<\sum_2^{\infty}
\frac{1}{k^2}=\frac{\pi^2}{6}-1=0.64\ldots<1,\] следовательно,
$|s_n(z)|>0$ при всех $|z|\le 1$ и $n\in \mathbb{N}$. Тем самым,
$n_0(f)=1$. Кроме того,
\[2-\pi^2/6\le |s_n(z)|\le \pi^2/6, \qquad |z|\le 1, \quad n\in \mathbb{N}.\]
Имеем также (см. (\ref{sum_k^(-s)})) \[|p(z)|\le |s_n(z)|\le
\frac{\pi^2}{6}, \qquad |f(z)|\le \sum_2^{\infty}
\frac{1}{k^2}=\frac{\pi^2}{6}-1, \qquad |z|\le 1;\]
\[|R_n(z)|\le \sum_{n+2}^{\infty}\frac{|z|^{k-1}}{k^2}\le
\frac{|z|^{n+1}}{n+1}, \qquad |z|\le 1, \qquad n\in \mathbb{N}.\]
Перейдём к оценкам производных. Очевидно,
\[|R_n'(z)|\le \sum_{k=n+1}^{\infty}\frac{k|z|^{k-1}}{(k+1)^2}\le
\sum_{k=n+1}^{\infty}\frac{|z|^{k-1}}{k+1}\le
\frac{1}{n+2}\cdot\frac{|z|^{n}}{1-|z|}, \qquad |z|<1, \qquad n\in
\mathbb{N}.\] Далее, $|p'(z)|\le n \pi^2/6$ при $|z|\le 1$ (см.
(\ref{|p'|})). Наконец, оценим $P(z)$ при $|z|\le 1$:
\[|P(z)|=|s_n(z)+z^m p(z)|\ge (1-|z|^m)|s_n(z)|\ge
(1-|z|^m)(2-\pi^2/6).\]

Применим полученные неравенства к оценке разности (\ref{P'/P-f}),
учитывая, что $2n+1\ge m\ge n+1$, в частности, $|z|^{m-1}\le
|z|^n$ при $|z|<1$:
\[\left|\frac{P'(z)}{P(z)}-f(z)\right|\le \frac{|z|^n}{(1-|z|^{n+1})(2-\pi^2/6)}
\left(\left(2n+1+\frac{\pi^2}{6}-1\right)\frac{\pi^2}{6}+n\frac{\pi^2}{6}+\right.\]
\[\left.+\frac{1}{n+2}\cdot\frac{1}{1-|z|}+
\left(\frac{\pi^2}{6}-1\right)\frac{1}{n+1}\right), \qquad
|z|<1.\] Поскольку $n\ge 1$, $(\pi^2/6)^2-1<11/6$ и $\pi^2<10$,
приходим к (\ref{result(|fj|<...)}):
\[\left|\frac{P'(z)}{P(z)}-f(z)\right|\le \frac{15|z|^n}{1-|z|^{n+1}}
\left(n+\frac{2}{1-|z|}\right), \qquad |z|<1.\] Теорема 1.3
доказана.

\end{document}